\documentclass[a4paper]{article}
\usepackage{amssymb}
\usepackage{mathrsfs}
\usepackage{float}
\usepackage{amsmath}
\usepackage{amsfonts,amssymb}
\usepackage{graphicx,color}
\usepackage{mathrsfs}
\catcode`\@=11 \@addtoreset{equation}{section}
\def\theequation{\thesection.\arabic{equation}}
\catcode`\@=12

\newtheorem{Theorem}{Theorem}[section]
\newtheorem{Definition}{Definition}[section]
\newtheorem{Lemma}{Lemma}[section]

\newtheorem{Corollary}{Corollary}[section]
\newtheorem{Remark}{Remark}[section]

\def\IRn{{I \hskip -1.6mm I\hskip -1.5mm R}}
\def\IR{{I \hskip -1.9mm I\hskip -1.7mm R}}
\def\2{{I \hskip -1.0mm I}}
\def\3{{I \hskip -1.0mm I\hskip -1.0mm I}}
\def\4{{I \hskip -0.9mm V}}
\def\6{{V \hskip -1.35mm I}}

\def\wt{\widetilde}

\def\ss{\scriptstyle}

\begin{document}
\title{}
\author{}
\date{}
\maketitle \catcode`\@=11 \@addtoreset{equation}{section}
\def\theequation{\thesection.\arabic{equation}}
\def\IRn{{I \hskip -1.6mm I\hskip -1.5mm R}}
\def\IR{{I \hskip -1.9mm I\hskip -1.7mm R}}
\def\wt{\widetilde}
\def\ss{\scriptstyle}
\vskip-3cm

\centerline{\Large\bf Lifespan of Classical Solutions to
Quasi-linear}

\bigbreak \centerline{\Large\bf Hyperbolic Systems with Small BV
Normal Initial Data}

\vskip 1.2cm \centerline{\large\bf Wen-Rong Dai \footnote{Center of
Mathematical Sciences, Zhejiang University, Hangzhou 310027, China,
E-mail: wrdai@126.com.}$\;$\footnote{This work was supported partly
by the National Science Foundation of China under Grant No. 10671124
and NO. 10771187.}}
\date{ }
\maketitle \vskip0.6cm

\centerline{\bf Abstract}
\begin{quote} In this paper, we first give a lower bound of the lifespan and some estimates of
classical solutions to the Cauchy problem for general quasi-linear
hyperbolic systems, whose characteristic fields are not weakly
linearly degenerate and the inhomogeneous terms satisfy Kong's
matching condition. After that, we investigate the lifespan of the
classical solution to the Cauchy problem and give a sharp limit
formula. In this paper, we only require that the initial data are
sufficiently small in the $L^1$ sense and the BV sense.
 \vskip0.5cm
{\it Key Words:}\ \  Quasi-linear hyperbolic system; Classical
solutions; Weakly linear degeneracy; Matching condition; Normalized
Coordinates; Blow-up; Lifespan. \vskip0.5cm {\it 2000 Mathematics
Subject Classification:} \ \ 35L45; 35L60; 35L40.
\end{quote}

\section{Introduction and main results}

\setcounter{equation}{0} Consider the following quasi-linear
hyperbolic system of first order
\begin{eqnarray} \frac {\partial u}{\partial t}+A(u)\frac {\partial u}{\partial
x}=B(u),
\end{eqnarray}
where $u=(u_1,\cdots,u_n)^{T}$ are the unknown vector-valued
functions of $(t,x)$, $A(u)=(a_{ij}(u))$ is an $n\times n$ matrix
and $B(u)=(B_1(u),B_2(u),\cdots,B_n(u))^T$ are $n-$dimensional
vector-valued functions.

By hyperbolicity, for any given $u$ on the domain under
consideration, $A(u)$ has $n$ real eigenvalues
$\lambda_1(u),\cdots,\lambda_n(u)$ and a complete system of left
(resp. right) eigenvectors $l_1(u),\cdots,l_n(u)\; $ $(resp.\;
r_1(u),\cdots,r_n(u))$. In this paper, we assume that (1.1) is a
strictly hyperbolic system, i.e.,
\begin{equation}
\lambda_1(u)<\lambda_2(u)<\cdots<\lambda_n(u).
\end{equation}
Without loss of generality, we suppose that on the domain under
consideration
\begin{eqnarray*}
l_{i}(u)r_j(u)\equiv \delta_{ij},\quad r_{i}^{T}(u)r_i(u)\equiv
1\;\; (i,j=1,\cdots,n),
\end{eqnarray*}
where $\delta_{ij}$ stands for the Kronecker's symbol.

The following definitions come from Kong \cite{k1}.
\begin{Definition}
The $i-th$ characteristic $\lambda_i(u)$ is {\bf weakly linearly
degenerate}, if, along the $i-$th characteristic trajectory
$u=u^{(i)}(s)$ passing through $u=0$, defined by
$\frac{du^{(i)}(s)}{ds}=r_i(u^{(i)}(s)),\;s=0:u=0,$ we have
\begin{equation}
\nabla\lambda_i(u)r_i(u)\equiv 0,\;\;\forall\;|s|\;\; small,
\end{equation}
namely
\begin{equation}
\lambda_i(u^{(i)}(s))\equiv \lambda(0),\;\;\forall\;|s|\;\; small.
\end{equation}
If all characteristics $\lambda_i(u)\;\;(i=1,2,\cdots,n)$ are weakly
linearly degenerate, then the system (1.1) is called {\bf weakly
linearly degenerate}.
\end{Definition}
\begin{Definition}
The inhomogeneous term $B(u)$ is called to be satisfied {\bf the
matching condition}, if, along the $i-$th characteristic trajectory
$u=u^{(i)}(s)$ passing through $u=0$, $B(u)\equiv 0$, i.e.,
\begin{equation}
B(u^{(i)}(s))\equiv 0,\;\;\forall \;|s| \;\;\;small.
\end{equation}
\end{Definition}
\begin{Definition}
If there exists a sufficiently smooth invertible transformation
$u=u(\tilde u)\;\;(u(0)=0)$ such that in the $\tilde u-$ space, for
each $i=1,2,\cdots,n$, the $i-$th characteristic trajectory passing
through $\tilde u=0$ coincides with the $\tilde u-$axis at least for
$|\tilde u_i|$ small, namely
\begin{equation}
\tilde r_i(\tilde u_ie_i)\equiv e_i,\;\;\forall\;|\tilde
u_i|\;\;small,
\end{equation}
where $e_i=(0,\cdots,0,\stackrel {(i)}{1},0,\cdots,0)^T$. Such a
transformation is called {\bf a normalized transformation} and the
corresponding unknown variables $\tilde u=(\tilde u_1,\tilde
u_2,\cdots,\tilde u_n)^T$ are called {\bf normalized variables} or
{\bf normalized coordinates}.
\end{Definition}

If the system (1.1) is strictly hyperbolic, then there always exists
the normalized transformation (cf. \cite{k2}). In this paper, for
the sake of simplicity, we assume that the unknown variables $u$ are
already normalized variables. That is to say,
\begin{equation}r_i(u_ie_i)\equiv e_i.
\end{equation}
It is easy to see
\begin{equation}r_i(0)=e_i,\;\;l_i(0)=e_i^T.\end{equation}
At the same time, (1.4) and (1.5) can be deduced to
\begin{equation}
\lambda_i(u_ie_i)\equiv 0
\end{equation}
and
\begin{equation}
B(u_ie_i)=0 \end{equation} respectively.

We consider the Cauchy problem of the hyperbolic system (1.1) with
the following initial data
\begin{eqnarray}
t=0:u(0,x)=f(\epsilon,x),
\end{eqnarray}
where $f(\epsilon,x)$ is a $C^1$ vector-valued function of
$\epsilon,\;x$ such that \begin{eqnarray} f(0,x)\equiv
0,\;\;\frac{\partial^2 f}{\partial\epsilon\partial
x}(\epsilon,\cdot)\in \left(C^r[0,\epsilon_1]\right)^n, \;0<r\le 1,
\end{eqnarray}
where $\epsilon\in [0,\epsilon_1]$, $\epsilon_1$ is a sufficiently
small positive constant. Then we know that
\begin{eqnarray}\underset{\epsilon\rightarrow
0^+}{\lim}\frac{f(\epsilon,x)}{\epsilon}=\frac{\partial f}{\partial
\epsilon}(0,x)\triangleq \psi(x)\in
\left(C^1(\mathbb{R})\right)^n.\end{eqnarray}

For the case that the initial data $f(\epsilon,x)$ satisfies the
following decay property: there exists a constant $\mu>0$ such that
\begin{eqnarray}
\varrho\stackrel{\triangle}{=}\sup_{x\in\mathbb{R}}\left\{(1+|x|)^{1+\mu}\left(|f(\epsilon,x)|+\left|\frac{\partial
f}{\partial x}(\epsilon,x)\right|\right)\right\}<+\infty
\end{eqnarray}
is sufficiently small, by means of the normalized coordinates Li et
al proved that the Cauchy problem (1.1) and (1.11) admits a unique
global classical solution, provided that the system (1.1) is weakly
linearly degenerate (see \cite{lk}-\cite{lzk2} and \cite{k2}). Kong
and Yang \cite{ky} studied the asymptotic behavior of the classical
solution. In their works, the condition $\mu
> 0$ is essential. If $\mu = 0$, a counterexample was constructed by
Kong \cite{k1} showing that the classical solution may blow up in a
finite time, even when the system (1.1) is weakly linearly
degenerate.

For the quasi-linear strictly hyperbolic system with linearly
degenerate characteristic fields, A. Bressan \cite{b} proved the
global existence of classical solution with initial data of small BV
norm. If the characteristic fields are weakly linearly degenerate,
Zhou \cite{z} proved the global existence of classical solution with
initial data of small $L^1$ norm and BV norm. Dai and Kong \cite{dk}
and Dai \cite{d1} studied the asymptotic behavior of the classical
solution.

When system (1.1) is not weakly linearly degenerate, there exists a
nonempty set $J\subseteq \{1,2,\cdots,n\}$ such that $\lambda_i(u)$
is not weakly linearly degenerate if and only if $i\in J$.

Noting (1.4), we observe that for any fixed $i\in J$, either there
exists an integer $\alpha_i\ge 0$ such that
\begin{eqnarray}\left.\frac {d^l\lambda_i\left(u^{(i)}(s)\right)}{ds^l}\right|_{s=0}=0\;\;
\left(l=1,\cdots,\alpha_i\right),\quad \mbox {but}\quad\left.\frac
{d^{\alpha_i+1}\lambda_i\left(u^{(i)}(s)\right)}{ds^{\alpha_i+1}}\right|_{s=0}\ne
0,\end{eqnarray} or
\begin{eqnarray}\left.\frac {d^l\lambda_i\left(u^{(i)}(s)\right)}{ds^l}\right|_{s=0}=0\quad
\left(l=1,2,\cdots\right).\end{eqnarray} In the case that (1.16)
holds, we define $\alpha_i=+\infty$.

For the normalized coordinates, conditions (1.15) and (1.16) simply
reduce to
\begin{equation*}
\frac{\partial ^l\lambda_i}{\partial
u_i^l}(0)=0\;\;(l=1,\cdots,\alpha_i),\;\;\mbox{but}\;\;\frac{\partial
^{\alpha_i+1}\lambda_i}{\partial u_i^{\alpha_i+1}}(0)\neq 0
\end{equation*}
and
\begin{equation*}
\frac{\partial ^l\lambda_i}{\partial u_i^l}(0)=0\;\;(l=1,2,\cdots)
\end{equation*}
respectively.

Our first goal in this paper is to give the following uniform {\it a
priori} estimates of the classical solution to the Cauchy problem
(1.1) and (1.11).
\begin{Theorem}
Suppose that the system (1.1) is strictly hyperbolic, $A\left
(u\right), B(u)$ is suitably smooth in a neighborhood of $u=0$ and
$B(u)$ satisfies the matching condition, suppose furthermore that
the initial data (1.11) satisfies
\begin{equation}
\int_{-\infty}^{+\infty}\left|\frac{\partial f}{\partial
x}(\epsilon,x)\right|dx\le
K_1\epsilon,\;\;\int_{-\infty}^{+\infty}|f(\epsilon,x)|dx\le
\frac{K_2}{M+1}\epsilon,
\end{equation}
where $\epsilon$ is a sufficiently small positive constant and
$K_1,K_2$ and $M\triangleq
\underset{x\in\mathbb{R}}{\sup}\left|\frac{\partial \psi}{\partial
x}(x)\right|$ are constants independent of $\epsilon$. Suppose
finally that system (1.1) is not weakly linearly degenerate and
\begin{eqnarray}\alpha =\min\left\{\alpha_i\mid i\in J\right\}<\infty ,\end{eqnarray}
where $\alpha_i$ is defined by (1.15)-(1.16). Then, on the existence
domain $[0,T]\times \mathbb{R}$ of the $C^1$ solution
$u=u\left(t,x\right)$, there exist positive constants
$K_3,\;K_4,\;K_5,\;K_6$ independent of $\epsilon,\; M,\;T$ such that
\begin{eqnarray}
V_1(T),\;\tilde V_1(T)\le
K_3(\epsilon+\epsilon^{\alpha+2}T),\;W_1(T),\; \tilde{W}_1(T),\;
U_{\infty}(T),\;V_{\infty}(T)\le K_3\epsilon,
\end{eqnarray}
where
\begin{equation}
T\epsilon^{\alpha+\frac{3}{2}}\le K_4
\end{equation}
and
\begin{equation}
W_{\infty} (T)\le K_5\epsilon,
\end{equation}
where
\begin{equation}
T\epsilon^{\alpha+1}\le K_6.
\end{equation}
In (1.19) and (1.21), $V_1(T),\;\tilde V_1(T),\;W_1(T),\;
\tilde{W}_1(T),\; U_{\infty}(T),\;V_{\infty}(T),\;W_{\infty} (T)$
are defined as follows: For any fixed $T\ge 0$,
\begin{eqnarray*}
U_{\infty}(T)=\sup_{0\le t\le T}\sup_{x\in\mathbb{R}}|u(t,x)|,\quad
V_{\infty}(T)=\sup_{0\le t\le T}\sup_{x\in\mathbb{R}}|v(t,x)|,
\end{eqnarray*}
\begin{eqnarray*}
W_{\infty}(T)=\sup_{0\le t\le T}\sup_{x\in\mathbb{R}}|w(t,x)|,
\end{eqnarray*}
\begin{eqnarray*}
V_1(T)=\sup_{0\le t\le T}\int_{-\infty}^{+\infty}|v(t,x)|dx,\quad
W_1(T)=\sup_{0\le t\le T}\int_{-\infty}^{+\infty}|w(t,x)|dx,
\end{eqnarray*}
\begin{eqnarray*}
\tilde V_1(T)=\max_{i\neq
j}\sup_{\tilde{C}_j}\int_{\tilde{C}_j}|v_i(t,x)|dt,\quad \tilde
W_1(T)=\max_{i\neq
j}\sup_{\tilde{C}_j}\int_{\tilde{C}_j}|w_i(t,x)|dt,
\end{eqnarray*}
where $|\cdot|$ stands for the Euclidean norm in $\mathbb{R}^n$,
$v=(v_1,\cdots,v_n)^T$ and $w=(w_1,\cdots,w_n)^T$ in which
$v_i=l_i(u)u$ and $w_i=l_i(u)u_x$ are defined by (2.1) in \S 2,
$\tilde{C}_j$ stands for any given $j$-th characteristic on the
domain $[0,T]\times\mathbb{R}$.
\end{Theorem}
\begin{Remark} By (1.21)-(1.22), we know that the life span of the classical solution $\tilde T(\epsilon)\ge K_6\epsilon^{-(\alpha+1)}$.
It is obvious that (1.14) implies (1.17). Therefore, Theorem 1.1 is
a generalization of corresponding results of Li et al \cite{lzk2}
and Kong \cite{k2} where the decay initial data was considered.
\end{Remark}

For the critical case, i.e., in (1.18), $\alpha=+\infty$, from
Theorem 1.1 and its proof in \S 3, we can easily get the following
corollary.
\begin{Corollary}
Assume that the assumptions except (1.18) in Theorem 1.1 hold. In
(1.18), we assume that $\alpha=+\infty$. Then, for any given integer
$N\ge 1$, there exists $\epsilon_0=\epsilon_0(N)>0$ so small that
for any fixed $\epsilon\in (0,\epsilon_0]$, the lifespan $\tilde
T(\epsilon)$ of the $C^1$ solution $u=u(t,x)$ to the Cauchy problem
(1.1) and (1.11) satisfies
$$\tilde T(\epsilon)\ge C_N\epsilon^{-N},$$
where $C_N$ is a positive constant independent of $\epsilon$.
\end{Corollary}

Next we consider the blow-up of the classical solution to the Cauchy
problem of the hyperbolic system (1.1) with the initial data (1.11).
If the hyperbolic system (1.1) is not weakly linearly degenerate, Li
et al \cite{lzk2} and Kong \cite{k2} estimated the lifespan of
classical solution to the Cauchy problem (1.1) with the special
initial data $u(0,x)=\epsilon \phi(x)$ which satisfies the following
decay property: there exists a constant $\mu>0$ such that
\begin{eqnarray}\varrho\stackrel{\triangle}{=}\sup_{x\in\mathbb{R}}\{(1+|x|)^{1+\mu}(|\phi(x)|+|\phi^{\prime}(x)|)\}<+\infty\end{eqnarray}
and the zero or matching inhomogeneous term $B(u)$.

Our second goal is to investigate the lifespan of classical solution
to the Cauchy problem (1.1) and (1.11) when the system (1.1) is not
weakly linearly degenerate.

\begin{Theorem} Suppose that the assumptions in Theorem 1.1 hold.
Let
\begin{eqnarray}J_1=\left\{i\,|\,\,i\in J,\,\alpha_i=\alpha\;\right\}\ne\emptyset.\end{eqnarray}
If there exists $i_0\in J_1$ and a point $x_0\in \mathbb{R}$ such
that
\begin{eqnarray}\frac{\partial^{\alpha+1}\lambda_{i_0}}{\partial u_{i_0}^{\alpha+1}}(0)\psi_{i_0}^{\alpha}\psi_{i_0}^{\prime}\left(x_0\right)<0,\end{eqnarray}
where $\psi(x)\in \left(C^1(\mathbb{R})\right)^n$ is defined in
(1.13), then there exists $\epsilon_0>0$ so small that for any fixed
$\epsilon\in\left (0,\epsilon_0\right]$, the first order derivatives
of the $C^1$ solution $u=u\left(t,x\right)$ to the Cauchy problem
(1.1) and (1.11) must blow up in a finite time and the lifespan
$\tilde {T}\left(\epsilon\right)$ of $u=u\left(t,x\right)$ satisfies
\begin{eqnarray}\lim_{\epsilon\rightarrow 0^+}\left(\epsilon^{\alpha+1} \tilde T
(\epsilon )\right)^{-1}=\max_{i\in J_1}\sup_{x\in
\mathbb{R}}\left(-\frac{1}{\alpha!}\frac{\partial^{\alpha+1}\lambda_i}{\partial
u_i^{\alpha+1}}(0)\psi_{i}^{\alpha}\left(x\right)\psi_i^{\prime}(x)\right).
\end{eqnarray}
\end{Theorem}
\begin{Remark} It is obvious that the decay property (1.23) implies (1.17). Therefore, Theorem
1.2 is a generalization of responding results of Li et al
\cite{lzk2} and Kong \cite{k2} and results of L. H\"ormander
\cite{h}, John \cite{j}, Liu \cite{liu} where the decay initial data
and the compactly supported initial data are considered
respectively.
\end{Remark}
\begin{Remark} For the hyperbolic system (1.1) with constant multiple
characteristic fields, we can obtain the similar results in Theorem
1.1 and Theorem 1.2 if we prove them as in this paper and in
\cite{d1}-\cite{d2}.
\end{Remark}
\begin{Remark} Similar to Kong and Li \cite{kl}, if  along $i-$th
characteristic $x=x_i(t,y)$, $w_i(t,x_i(t,y))=l_i(u)u_
x(t,x_i(t,y))$ blow up at the lifespan $\tilde T(\epsilon)$, then we
have $$w_i(t,x_i(t,y))=O((\tilde
T(\epsilon)-t)^{-1}),\quad\mbox{when}\;\;t\rightarrow\tilde
T(\epsilon)^-.$$\end{Remark}
\begin{Remark} For the conservation laws, shock will appear (see Kong \cite{k3}).
\end{Remark}

This paper is organized as follows. In \S 2, we recall John's
formula on the decomposition of waves with some supplements for the
hyperbolic system (1.1). Then we give some uniform {\it a priori}
estimates for the Cauchy problem (1.1) and (1.11) and prove Theorem
1.1 in \S 3. In \S 4, we obtain some important uniform estimates by
making use of an invertible characteristics' transformation of the
hyperbolic system (1.1). Finally, we investigate the lifespan of the
classical solution to the Cauchy problem (1.1) and (1.11) and give
the proof of Theorem 1.2 in \S 5.

\section{Preliminaries and Decomposed Formulas of Waves}

For the sake of completeness, in this section we briefly recall
John's formula on the decomposition of waves with some supplements
for the hyperbolic system (1.1), which play an important role in our
proof.

Let
\begin{equation} v_i=l_i(u)u,\quad
w_i=l_i(u)u_x\quad (i=1,\cdots ,n)\end{equation} and
\begin{eqnarray}
b_i(u)=l_i(u)B(u)\quad (i=1,2,\cdots,n).
\end{eqnarray}
Then we have
\begin{equation}u=\sum^n_{k=1}v_kr_k(u),\quad
u_x=\sum^n_{k=1}w_kr_k(u)\end{equation} and
\begin{eqnarray}
B(u)=\sum_{k=1}^{n}b_k(u)r_k(u).
\end{eqnarray}
Let
\begin{equation}\frac d{d_it}=\frac {\partial}{\partial
t}+\lambda_i(u)\frac {\partial}{
\partial x}\end{equation}
be the directional derivative along the $i$-th characteristic. We
have (see \cite{lkz}-\cite{lzk2} or \cite{k2})
\begin{eqnarray}
\frac
{dv_i}{d_it}=\sum_{j,k=1}^{n}\beta_{ijk}(u)v_jw_k+\sum_{j,k=1}^{n}\nu_{ijk}(u)v_jb_k(u)+b_i(u)\stackrel
{\triangle}{=} F_i(t,x)
\end{eqnarray}
and
\begin{eqnarray}
\frac
{dw_i}{d_it}=\sum_{j,k=1}^{n}\gamma_{ijk}(u)w_jw_k+\sum_{j,k=1}^{n}\sigma_{ijk}(u)w_jb_k(u)+(b_i(u))_x\stackrel
{\triangle}{=} G_i(t,x),
\end{eqnarray}
where
\begin{eqnarray}
\beta_{ijk}(u)&=&(\lambda_k(u)-\lambda_i(u))l_i(u)\nabla
r_j(u)r_k(u),\\
\nu_{ijk}(u)&=&-l_i(u)\nabla r_j(u)r_k(u),\\
\gamma_{ijk}(u)&=&(\lambda_k(u)-\lambda_j(u))l_i(u)\nabla
r_j(u)r_k(u)-\nabla \lambda_j(u)r_k(u)\delta_{ij},\\
\sigma_{ijk}(u)&=&l_i(u))(\nabla r_k(u)r_j(u)-\nabla r_j(u)r_k(u)).
\end{eqnarray} Equivalently we also get
\begin{eqnarray}
\frac{\partial v_i}{\partial t}+\frac{\partial
(\lambda_i(u)v_i)}{\partial x}&=&\sum_{j,k=1}^{n}\tilde
{\beta}_{ijk}(u)v_jw_k+\sum_{j,k=1}^{n}\nu_{ijk}(u)v_jb_k(u)+b_i(u)\nonumber\\&\triangleq&
\tilde {F}_i(t,x),\\
d[v_i(dx-\lambda_i(u)dt)]&=&\left[\sum_{j,k=1}^{n}\tilde
{\beta}_{ijk}(u)v_jw_k+\sum_{j,k=1}^{n}\nu_{ijk}(u)v_jb_k(u)+b_i(u)\right]dt\wedge
dx  \nonumber\\ &\stackrel {\triangle}{=} &\tilde {F}_i(t,x)
dt\wedge dx
\end{eqnarray} and
\begin{eqnarray}
\frac{\partial w_i}{\partial t}+\frac{\partial
(\lambda_i(u)w_i)}{\partial x}&=&\sum_{j,k=1}^{n}\tilde
{\gamma}_{ijk}(u)w_jw_k+\sum_{j,k=1}^{n}\sigma_{ijk}(u)w_jb_k(u)+(b_i(u))_x\nonumber\\
&\stackrel {\triangle}{=} &\tilde {G}_i(t,x),\\
d[w_i(dx-\lambda_i(u)dt)]&=&\left[\sum_{j,k=1}^{n}\tilde
{\gamma}_{ijk}(u)w_jw_k+\sum_{j,k=1}^{n}\sigma_{ijk}(u)w_jb_k(u)+(b_i(u))_x\right]dt\wedge
dx \nonumber\\ &\stackrel {\triangle}{=} &\tilde {G}_i(t,x) dt\wedge
dx,
\end{eqnarray}
where
\begin{eqnarray}
\tilde {\beta}_{ijk}(u)&=&\beta_{ijk}(u)+\nabla
\lambda_i(u)r_k(u)\delta_{ij},\\
\tilde {\gamma}_{ijk}(u)&=&\gamma_{ijk}(u)+\frac {1}{2}[\nabla
\lambda_j(u)r_k(u)\delta_{ij}+\nabla \lambda_k(u)r_j(u)\delta_{ik}].
\end{eqnarray}

From (2.8), (2.10) and (2.16)-(2.17), we see that
\begin{eqnarray}
\beta_{iji}(u)&\equiv& 0,\;\;\tilde {\gamma}_{ijj}(u)\equiv 0,\quad
\forall\; i,j\in \{1,2,\cdots,n\},\quad \forall\;|u|\;\; small,\\
\gamma_{ijj}(u)&\equiv& 0,\;\;\tilde {\beta}_{iji}(u)\equiv 0,\quad
\forall\; j\ne i,\quad \forall\;|u|\;\; small.
\end{eqnarray}

As we already assume that $u$ are the normalized coordinates, making
use of (1.7), the following relations hold (see \cite {k2}):
\begin{eqnarray}
\beta_{ijj}(u_je_j)&\equiv& 0,\quad\nu_{ijj}(u_je_j)\equiv
0,\quad\sigma_{ijj}(u_je_j) \equiv 0,\;\; |u_j| \;\hbox
{small},\; \forall\; i,j,\\
\tilde {\beta}_{ijj}(u_je_j)&\equiv&
0,\quad\quad\quad\quad\quad\quad\forall\; |u_j| \;\hbox {small},\;
\forall\; i\neq j.
\end{eqnarray}

When the inhomogeneous term $B(u)$ satisfies the matching condition,
then in the normalized coordinates $u$ (see \cite{k2}),
\begin{eqnarray}
b_i(u)&=&\sum_{j\neq k}b_{ijk}(u)u_ju_k, \quad
 \forall\; |u| \;\hbox {small},\; \forall\; i\in \{1,2,\cdots,n\},\\
 (b_i(u))_x&=&\sum_{k=1}^{n}\tilde {b}_{ik}(u)w_k,
\end{eqnarray}
where $b_{ijk}(u)$ is a $C^1$ function and $\tilde {b}_{ik}(u)=
\overset{n}{\underset{l=1}{\sum}}\frac {\partial b_i(u)}{\partial
u_l}r_{kl}(u)$ satisfies that
\begin{eqnarray}
\tilde {b}_{ik}(u_ke_k)\equiv 0,\quad \forall\; |u_k|\; \hbox
{small}, \; \forall\; k\in\{1,\cdots,n\}.
\end{eqnarray}

\section{Uniform Estimates---Proof of Theorem 1.1}

In this section, we shall establish some uniform estimates under the
assumptions in Theorem 1.1 and give the proof of Theorem 1.1.

First we recall some basic $L^1$ estimates. They are essentially due
to Schartzman \cite{s1}, \cite{s2} and Zhou \cite{z}.

\begin{Lemma}
Let $\phi=\phi(t,x)\in C^1$ satisfy
\begin{equation*}
\phi_t+(\lambda(t,x)\phi)_x=F(t,x),\;\;0\le t\le
T,x\in\mathbb{R},\;\;\phi(0,x)=g(x),
\end{equation*}
where $\lambda\in C^1$. Then
\begin{equation*}
\int_{-\infty}^{+\infty}|\phi(t,x)|dx\le\;\int_{-\infty}^{+\infty}|g(x)|dx+\int_0^T\int_{-\infty}^{+\infty}|F(s,x)|dsdx,\;\;\forall\;t\le
T,
\end{equation*}
provided that the right hand side of the inequality is bounded.
\end{Lemma}

\begin{Lemma}
Let $\phi=\phi(t,x)$ and $\psi=\psi(t,x)$ be $C^1$ functions
satisfying
\begin{equation*}
\phi_t+(\lambda(t,x)\phi)_x=F(t,x),\;\;0\le t\le
T,x\in\mathbb{R},\;\;\phi(0,x)=g_1(x),
\end{equation*}
and
\begin{equation*} \psi_t+(\mu(t,x)\phi)_x=G(t,x),\;\;0\le t\le
T,x\in\mathbb{R},\;\;\phi(0,x)=g_2(x),
\end{equation*}
respectively, where $\lambda,\;\mu\in C^1$ such that there exists a
positive constants $\delta_0$ independent of $T$ verifying
\begin{equation*}
\mu(t,x)-\lambda(t,x)\ge \delta_0,\;\;0\le t\le
T,\;\;x\in\mathbb{R}.
\end{equation*}
Then \begin{eqnarray*}\begin{array}{lll}
\int_0^T\int_{-\infty}^{+\infty}|\phi(t,x)||\psi(t,x)|dxdt &\le
&C\left(\int_{-\infty}^{+\infty}|g_1(x)|dx+\int_0^T\int_{-\infty}^{+\infty}|F(t,x)|dxdt\right)\times\vspace{2mm}\\
&&\;\;\;\;\left(\int_{-\infty}^{+\infty}|g_2(x)|dx+\int_0^T\int_{-\infty}^{+\infty}|G(t,x)|dxdt\right),
\end{array}\end{eqnarray*}
provided that the two factors on the right hand side of the
inequality is bounded.
\end{Lemma}

By the existence and uniqueness of local $C^1$ solution to the
Cauchy problem, in order to prove Theorem 1.1, it suffices to
establish {\it a prior} estimates on the $C^{0}$ norm of $u$ and
$\frac {\partial u}{\partial x}$ on the existence domain of $C^1$
solution $u=u(t,x)$.

By (1.2), there exist positive constants $\delta_0$, $\delta_1$ and
$\delta$ such that
\begin{equation}|\lambda_{i}(u)-\lambda_j(v)|\geq \delta_0,\;\;|\lambda_i\left(u\right)-\lambda_i\left(v\right)|\le \delta_1,
\quad\forall\;|u|,\,|v|\leq\delta,\quad\forall\; i\neq
j.\end{equation}

For the time being it is supposed that on the existence domain
$[0,T]\times \mathbb{R}$ of the $C^1$ solution $u=u\left(t,x\right)$
we have \begin{eqnarray} |u(t,x)|\le K_7\epsilon,\end{eqnarray}
where $K_7$ is a positive constant independent of
$\epsilon,\;t,\;x$. At the end of the proof of Theorem 1.1, we shall
explain that this hypothesis is reasonable. Then, (3.1) hold if we
take $\delta=K_7\epsilon$.

Introduce
\begin{eqnarray*}
Q_W(T)&=&\sum_{j\neq k}\int_{0}^{T}\int_{\mathbb{R}}|w_j(t,x)||w_k(t,x)|dtdx,\nonumber\\
Q_{VW}(T)&=&\sum_{j\neq
k}\int_{0}^{T}\int_{\mathbb{R}}|v_j(t,x)||w_k(t,x)|dtdx,\nonumber\\
Q_V(T)&=&\sum_{j\neq
k}\int_{0}^{T}\int_{\mathbb{R}}|v_j(t,x)||v_k(t,x)|dtdx.\nonumber\\
\end{eqnarray*}
As we already assume that $u$ are the normalized coordinates, by
(1.7) it can be easily seen that
\begin{equation*}
\sum_{i\neq j}|u_i|\le C_{1}\sum_{i\neq
j}|v_i|,\;\mbox{for\;fixed}\;\; j;\; \sum_{i\neq j}|u_iw_j|\le
C_1\sum_{i\neq j}|v_iw_j|;\;\;\;\sum_{i\neq j}|u_iu_j|\le
C_1\sum_{i\neq j}|v_iv_j|.
\end{equation*}
Here and hereafter $C_j\;(j=1,2,\cdots)$ stand for some positive
constants independent of $\epsilon,\;M,\;T$.

It follows from (2.12) and (2.18)-(2.24) that
\begin{eqnarray*}
\tilde F_i(t,x)&=&\sum_{j,k=1}^{n}\tilde
{\beta}_{ijk}(u)v_jw_k+\sum_{j,k=1}^{n}\nu_{ijk}(u)v_jb_k(u)+b_i(u)\\
&=&\sum_{j\neq k}\tilde {\beta}_{ijk}(u)v_jw_k+\sum_{j=1}^n(\tilde
{\beta}_{ijj}(u)-\tilde {\beta}_{ijj}(u_je_j))v_jw_j+\tilde
{\beta}_{iii}(u_ie_i)v_iw_i\\
&&+\sum_{j,k=1}^{n}\nu_{ijk}(u)v_j\sum_{p\neq
q}b_{kpq}(u)u_pu_q+\sum_{p\neq q}b_{ipq}(u)u_pu_q.
\end{eqnarray*}
On the other hand,
\begin{eqnarray*}
\tilde {\beta}_{iii}(u_ie_i)=\nabla
\lambda_i(u_ie_i)r_i(u_ie_i)=\frac{\partial \lambda_i}{\partial
u_i}(u_ie_i).
\end{eqnarray*}
Therefore, we have
\begin{equation}
|\tilde F_i(t,x)|\le C_2\left[\sum_{j\neq k}|v_jw_k|+\sum_{j\neq
k}|v_jv_k|+|u_i|^{\alpha}|v_iw_i|\right].
\end{equation}
By (2.14) and (2.18)-(2.24), we have
\begin{eqnarray*}
\tilde G_i(t,x)&=&\sum_{j,k=1}^{n}\tilde
{\gamma}_{ijk}(u)w_jw_k+\sum_{j,k=1}^{n}\sigma_{ijk}(u)w_jb_k(u)+(b_i(u))_x\\
&=&\sum_{j\neq k}\tilde
{\gamma}_{ijk}(u)w_jw_k+\sum_{j,k=1}^{n}\sigma_{ijk}(u)w_j\sum_{p\neq
q}b_{kpq}(u)u_pu_q\\
&&+\sum_{k=1}^n(\tilde b_{ik}(u)-\tilde b_{ik}(u_ke_k))w_k.
\end{eqnarray*}
Then we get
\begin{equation}
|\tilde G_i(t,x)|\le C_3\left[\sum_{j\neq k}|v_jw_k|+\sum_{j\neq
k}|w_jw_k|\right].
\end{equation}
By (2.12), (2.14), (3.3)-(3.4), it follows from Lemma 3.2 that
\begin{eqnarray}
Q_W(T)&\le & C_4\left(W_1(0)+\int_0^T\int_{\mathbb{R}}|\tilde
G(t,x)|dtdx\right)^2\nonumber\\&\le&C_4(W_1(0)+Q_W(T)+Q_{VW}(T))^2\nonumber\\
&\le&C_4(\epsilon+Q_W(T)+Q_{VW}(T))^2,\\
Q_V(T)&\le & C_4\left(V_1(0)+\int_0^T\int_{\mathbb{R}}|\tilde
F(t,x)|dtdx\right)^2\nonumber\\&\le&C_4(V_1(0)+Q_V(T)+Q_{VW}(T)+|U_{\infty}(T)|^{\alpha+1}W_1(T)\cdot T)^2\nonumber\\
&\le&C_4\left(\frac{\epsilon}{M+1}+Q_V(T)+Q_{VW}(T)+|U_{\infty}(T)|^{\alpha+1}W_1(T)\cdot T\right)^2,\\
Q_{VW}(T)&\le & C_4\left(V_1(0)+\int_0^T\int_{\mathbb{R}}|\tilde
F(t,x)|dtdx\right)\left(W_1(0)+\int_0^T\int_{\mathbb{R}}|\tilde
G(t,x)|dtdx\right)\nonumber\\&\le&C_4(V_1(0)+Q_V(T)+Q_{VW}(T)+|U_{\infty}(T)|^{\alpha+1}
W_1(T)\cdot T)
\nonumber\\&&\cdot(W_1(0)+Q_W(T)+Q_{VW}(T))\nonumber\\
&\le&C_4\left(\frac{\epsilon}{M+1}+Q_V(T)+Q_{VW}(T)+|U_{\infty}(T)|^{\alpha+1}
W_1(T)\cdot T\right)\nonumber\\&&\cdot(\epsilon+Q_W(T)+Q_{VW}(T)),
\end{eqnarray}
where $\tilde F=(\tilde F_1,\tilde F_2,\cdots,\tilde F_n)^T,\;\tilde
G=(\tilde G_1,\tilde G_2,\cdots,\tilde G_n)^T$.

We assume that the $j-$th characteristic $\tilde C_j$ intersects
$t=0$ with point $A$, intersects $t=T$ with point $B$. We draw an
$i-$th characteristic $\tilde C_i$ from $B$ downward and intersects
$t=0$ with point $C$. We rewrite (2.15) as
\begin{equation*}
d(|w_i(t,x)|(dx-\lambda_i(u)dt))=sgn(w_i)\tilde G_idtdx,\;\;a.e.
\end{equation*}
and integrate it in the region $ABC$ to get
\begin{eqnarray*}
\left|\int_{\tilde
C_j}|w_i(t,x)|(\lambda_j(u)-\lambda_i(u))dt\right|&\le&\int_A^C|w_i(0,x)|dx+\int\int_{ABC}|\tilde
G_i|dtdx \\ &\le& C_5(W_1(0)+Q_W(T)+Q_{VW}(T)).
\end{eqnarray*}
Noting (3.1), it follows that
\begin{equation*}
\int_{\tilde C_j}|w_i(t,x)|dt\le C_6(W_1(0)+Q_W(T)+Q_{VW}(T))\le
C_7(\epsilon+Q_W(T)+Q_{VW}(T)),\end{equation*} hence
\begin{equation}
\tilde W_1(T)\le C_8(\epsilon+Q_W(T)+Q_{VW}(T)).
\end{equation}
In a similar way, we can deduce from (2.13) that
\begin{equation}
\tilde V_1(T)\le
C_9\left[\frac{\epsilon}{M+1}+Q_V(T)+Q_{VW}(T)+|U_{\infty}(T)|^{\alpha+1}
W_1(T)\cdot T\right].
\end{equation}

It follows from (2.12) and Lemma 3.1 that
\begin{eqnarray*}
\int_{-\infty}^{+\infty}|v_i(T,x)|dx&\le&\int_{-\infty}^{+\infty}|v_i(0,x)|dx
+\int_0^T\int_{-\infty}^{+\infty}|\tilde
F_i(t,x)|dtdx\nonumber\\&\le&C_{10}V_1(0)+C_2\int_0^T\int_{-\infty}^{+\infty}\left[\sum_{j\neq
k}|v_jw_k|+\sum_{j\neq
k}|v_jv_k|+|u_i|^{\alpha}|v_iw_i|\right]dtdx\nonumber\\&\le&C_{11}\left[\frac{\epsilon}{M+1}+Q_V(T)+Q_{VW}(T)+|U_{\infty}(T)|^{\alpha+1}
W_1(T)\cdot T\right].
\end{eqnarray*}
That is to say,
\begin{eqnarray}
V_1(T)\le
C_{11}\left[\frac{\epsilon}{M+1}+Q_V(T)+Q_{VW}(T)+|U_{\infty}(T)|^{\alpha+1}
W_1(T)\cdot T\right].
\end{eqnarray}
In a similar way, it follows from (2.14) and Lemma 3.1 that
\begin{eqnarray}
W_1(T)\le C_{12}\left[\epsilon+Q_W(T)+Q_{VW}(T)\right].
\end{eqnarray}

It can be easily seen that
\begin{eqnarray}
U_{\infty}(T),\;\;V_{\infty}(T)\le C_{13}\underset{0\le t\le
T}{\sup}\int_{-\infty}^{+\infty}|u_x(t,x)|dx\le C_{14}W_1(T).
\end{eqnarray}

Thus, in order to prove (1.19) it suffices to show that we can
choose some constants $d_i\;(i=1,2,3,4,5)$ in such a way that for
any fixed $T_0\;(0\le T_0\le T)$ with
$T_0\epsilon^{\alpha+\frac{3}{2}}\le K_4$ such that
\begin{eqnarray}
&&V_1(T_0),\;\tilde V_1(T_0)\le
2d_1\epsilon+2d_2\epsilon^{\alpha+2}T_0,\;W_1(T_0)\le
2d_3\epsilon,\; \tilde{W}_1(T_0)\le 2d_4\epsilon,\nonumber\\
&&U_{\infty}(T_0),\;V_{\infty}(T_0)\le 2d_5\epsilon,
\end{eqnarray}
we have
\begin{eqnarray}
&&V_1(T_0),\;\tilde V_1(T_0)\le
d_1\epsilon+d_2\epsilon^{\alpha+2}T_0,\;W_1(T_0)\le d_3\epsilon,\;
\tilde{W}_1(T_0)\le d_4\epsilon,\nonumber\\ &&
U_{\infty}(T_0),\;V_{\infty}(T_0)\le d_5\epsilon.
\end{eqnarray}

Substituting (3.13) into (3.5)-(3.7), we have
\begin{eqnarray*}
Q_W(T_0)&\le&C_4(\epsilon+Q_W(T_0)+Q_{VW}(T_0))^2,\\
Q_V(T_0)&\le&C_4\left[\frac{\epsilon}{M+1}+Q_V(T_0)+Q_{VW}(T_0)+(2d_5)^{\alpha+1}(2d_3)K_4\epsilon^{\frac{1}{2}}\right]^2,\\
Q_{VW}(T_0)&\le&C_4\left[\frac{\epsilon}{M+1}+Q_V(T_0)+Q_{VW}(T_0)+(2d_5)^{\alpha+1}(2d_3)K_4\epsilon^{\frac{1}{2}}\right](\epsilon+Q_W(T_0)+Q_{VW}(T_0)).
\end{eqnarray*}
Denote $a_1=(2d_5)^{\alpha+1}(2d_3)K_4$. It follows that
\begin{eqnarray}
Q_W(T_0)\le C_4(1+3C_4a_1^2)^2\epsilon^2,\;Q_V(T_0)\le
2C_4a_1^2\epsilon,\;Q_{VW}(T_0)\le 2C_4a_1\epsilon^{\frac{3}{2}},
\end{eqnarray}
provided that $\epsilon$ is sufficiently small.

Furthermore, by making use of (3.15), from (3.8)-(3.12), we get
\begin{eqnarray*}
\tilde W_1(T_0)&\le& 2C_8\epsilon,\;\tilde V_1(T_0)\le
2C_9\left[(\frac{1}{M+1}+2C_4a_1^2)\epsilon+(2d_5)^{\alpha+1}(2d_3)\epsilon^{\alpha+2}T_0\right],\\
W_1(T_0)&\le&
2C_{12}\epsilon,\;V_1(T_0)\le2C_{11}\left[(\frac{1}{M+1}+2C_4a_1^2)\epsilon+(2d_5)^{\alpha+1}(2d_3)\epsilon^{\alpha+2}T_0\right],\;\\
U_{\infty}(T_0),V_{\infty}(T_0)&\le& 2C_{12}C_{14}\epsilon.
\end{eqnarray*}

If we take
\begin{eqnarray*}
d_3\ge 2C_{12},\;d_4\ge 2C_8,\;d_5\ge 2C_{12}C_{14}
\end{eqnarray*}
and
\begin{eqnarray*}
d_1\ge
2\max\{C_9,C_{11}\}\left[\frac{1}{M+1}+2C_4a_1^2\right],\;d_2\ge
2\max\{C_9,C_{11}\}(2d_5)^{\alpha+1}(2d_3),
\end{eqnarray*}
then we obtain (3.14). Thus, if we take
$K_3=\max\;\{d_1,d_2,d_3,d_4,d_5\}$, we obtain (1.19).

It follow from (2.7) that
\begin{eqnarray*}
w_i(t,x_i(t,y))=w_i(0,y)+\int_{\tilde C_i}G_i(t,x_i(t,y))dt,
\end{eqnarray*}
where $\tilde C_i$ is the $i-$th characteristic defined by
$$\frac{d
x_i(t,y)}{dt}=\lambda_i(u(t,x_i(t,y))),\;\;t=0:x_i(0,y)=y.$$ By
(2.7) and (2.18)-(2.24), we have
\begin{eqnarray*}
G_i(t,x)&=&\sum_{j,k=1}^{n}
\gamma_{ijk}(u)w_jw_k+\sum_{j,k=1}^{n}\sigma_{ijk}(u)w_jb_k(u)+(b_i(u))_x\\
&=&\sum_{j\neq
k}\gamma_{ijk}(u)w_jw_k+(\gamma_{iii}(u)-\gamma_{iii}(u_ie_i))w_i^2+\gamma_{iii}(u_ie_i)w_i^2\nonumber\\
&&+\sum_{j,k=1}^{n}\sigma_{ijk}(u)w_j\sum_{p\neq
q}b_{kpq}(u)u_pu_q+\sum_{k=1}^n(\tilde b_{ik}(u)-\tilde
b_{ik}(u_ke_k))w_k.
\end{eqnarray*}
On the other hand,
\begin{eqnarray*}
\gamma_{iii}(u_ie_i)=-\nabla
\lambda_i(u_ie_i)r_i(u_ie_i)=-\frac{\partial \lambda_i}{\partial
u_i}(u_ie_i).
\end{eqnarray*}
Therefore, we get
\begin{equation}
|G_i(t,x)|\le C_{15}\left[\sum_{j\neq
k}\left(|w_jw_k|+|v_jw_k|\right)+\sum_{j\neq
i}|v_jw_i^2|+|u_i|^{\alpha}|w_i|^2\right].
\end{equation}
Then we obtain
\begin{eqnarray}
W_{\infty}(T)&\le& C_{16}\left[W_{\infty}(0)+W_{\infty}(T)\tilde
W_1(T)+V_{\infty}(T)\tilde W_1(T)+W_{\infty}(T)\tilde
V_1(T)\right.\nonumber\\
&&\left.\quad\;\;\;\;+W_{\infty}(T)^2\tilde
V_1(T)+U_{\infty}(T)^{\alpha}W_{\infty}(T)^2\cdot
T\right]\nonumber\\
&\le& C_{17}\left[\epsilon+W_{\infty}(T)\tilde
W_1(T)+V_{\infty}(T)\tilde W_1(T)+W_{\infty}(T)\tilde
V_1(T)\right.\nonumber\\
&&\left.\quad\;\;\;\;+W_{\infty}(T)^2\tilde
V_1(T)+U_{\infty}(T)^{\alpha}W_{\infty}(T)^2\cdot T\right].
\end{eqnarray}

Thus, in order to prove (1.21) it suffices to show that we can
choose some constant $d_6$ in such a way that, for any fixed
$T_1\;(0\le T_1\le T)$ with $T_1\epsilon^{\alpha+1}\le K_6$,
\begin{eqnarray}
W_{\infty}(T_1)\le 2d_6\epsilon,
\end{eqnarray}
we have
\begin{eqnarray}
W_{\infty}(T_1)\le d_6\epsilon.
\end{eqnarray}
Substituting (1.19) and (3.18) into (3.17), we have
\begin{eqnarray*}
W_{\infty}(T_1)\le 2C_{17}[1+K_3^{\alpha}(2d_6)^2K_6]\epsilon.
\end{eqnarray*}
Hence, if $d_6\ge 4C_{17},\;K_6=\frac{1}{K_3^{\alpha}(2d_6)^2}$,
then we have (3.19). Therefore (1.21) is proved.

It follows from (1.19) that $U_{\infty}(T)\le K_7\epsilon$ where $T$
satisfies (1.20), provided that $\epsilon$ is sufficient small and
$K_7\ge K_3$. Then the hypothesis (3.2) is reasonable.  This proves
Theorem 1.1. $\quad\Box$

\section{Some important uniform estimates on classical solutions}

On the domain where the classical solution $u=u(t,x)$ of the Cauchy
problem (1.1) and (1.11) exists, we denote the $i$-th characteristic
passing through the point ($0,y$) by $x=\phi^{(i)} (t,y)$, which is
defined by
\begin{eqnarray}\frac {\partial \phi^{(i)} (t,y)}{\partial t}=\lambda_i\left(
u\left(t,\phi^{(i)} (t,y) \right)\right),\quad \phi^{(i)} (0,y)=y.
\end{eqnarray} Let
\begin{eqnarray}z^{(i)}(t,y)=u\left(t,\phi^{(i)} (t,y)\right).\end{eqnarray}
For the sake of simplicity, we omit the upper index $(i)$ of
$z^{(i)},\;\phi^{(i)}$ etc. in this section. Then from (1.1) we
easily have
\begin{eqnarray}
l_i(z)\partial_t z=b_i(z) \end{eqnarray} and
\begin{eqnarray}l_j(z)\partial_y
z=\frac{b_j(z)-l_j(z)\partial_t
z}{\lambda_j(z)-\lambda_i(z)}(\partial_y \phi),\quad \forall\; j\neq
i.
\end{eqnarray}
\begin{Theorem} Under the assumptions of Theorem
1.1, we know that $(\phi,z)=(\phi (t,y),z(t,y))$ is $C^1$ smooth
with respect to $(t,y)$ on the domain
\begin{eqnarray*}D(M_1)=\left\{(t,y)\mid 0\le t< \min\left\{\tilde
T(\epsilon),\;M_1\epsilon^{-(\alpha+1)}\right\},\;-\infty<y<\infty\right\},\end{eqnarray*}
provided that $\epsilon$ is sufficiently small, where $M_1$ is any
positive constant independent of $\epsilon,\;t,\;y$ and $\tilde
T(\epsilon)$ is the lifespan of the $C^1$ classical solution
$u=u(t,x)$ to the Cauchy problem (1.1) and (1.11). Moreover, we have
$\phi_{ty}\in C^0$ and the following estimates hold in the domain
$D(M_1)$:
\begin{eqnarray}
|\phi_t(t,y)|&\le& C_{18},\;|\phi_y(t,y)|\le
C_{18},\;|\phi_{ty}(t,y)|\le C_{18}\epsilon,\nonumber\\
|z(t,y)|&\le& C_{18}\epsilon,\;|z_t(t,y)|\le
C_{18}\epsilon,\;|z_y(t,y)|\le C_{18}\epsilon.\end{eqnarray} In
addition, in the domain $D(M_1)$,
\begin{eqnarray}\bar
w_j(t,y)\triangleq \frac{b_j(z)-l_j(z)\frac{\partial z}{\partial
t}}{\lambda_j(z)-\lambda_i(z)}\in C^0,\;\;|\bar w_j(t,y)|\le
C_{19}\epsilon,\;\;j\ne i.\end{eqnarray}
\end{Theorem}
\begin{Remark}
In the existence domain of the $C^1$ solution $u=u(t,x)$ to the
Cauchy problem (1.1) and (1.11), i.e., in the domain
$\left[0,\min\left\{\tilde
T(\epsilon),M_1\epsilon^{-(\alpha+1)}\right\}\right)$ $\times
(-\infty,+\infty)$, from (1.1) and (4.4) we have, along the $i-th$
characteristic $x=\phi(t,y)$ passing the point $(0,y)$,
\begin{eqnarray*}u(t,\phi(t,y))\equiv z(t,y)\end{eqnarray*}
and
\begin{eqnarray*}
\bar w_j(t,y)&=&\frac{b_j(u)-l_j(u)\left(u_t+\lambda_i(u)u_x\right)}{\lambda_j(u)-\lambda_i(u)}(t,\phi(t,y))\\
&=&\frac{b_j(u)-l_j(u)\left(-A(u)u_x+B(u)+\lambda_i(u)u_x\right)}{\lambda_j(u)-\lambda_i(u)}(t,\phi(t,y))\\
&=&w_j(t,\phi(t,y)),\;\;\quad\quad\forall\;j\ne i,
\end{eqnarray*}
where $u=u(t,x)$ is the $C^1$ smooth solution to the Cauchy problem
(1.1) and (1.11) and $w_j=l_j(u)u_x$ is defined by (2.1). It follows
from (4.6) that
\begin{eqnarray}|w_j(t,\phi(t,y))|\le C_{19}\epsilon,\;j\ne i,\;\;\mbox{if}\;\;t\in \left[0,\min\left\{\tilde
T(\epsilon),M_1\epsilon^{-1}\right\}\right),\;\;y\in\mathbb{R}.\end{eqnarray}
\end{Remark}

\noindent{\bf Proof.} It follows from (4.1) that
\begin{eqnarray*}\phi_{ty}(t,y)=\left(\lambda_i(u(t,\phi(t,y)))\right)_y=\sum_{j=1}^{n}\frac{\partial
\lambda_i(u)}{\partial u_j}(
u_j)_x\phi_y(t,y),\;\;\phi_y(0,y)=1.\end{eqnarray*} Then, we get
\begin{eqnarray*}\ln\left|\phi_y(t,y)\right|=\int_0^t\sum_{j=1}^{n}\frac{\partial
\lambda_i(u)}{\partial u_j}(u_j)_x(t,\phi(t,y))\end{eqnarray*}
Before the blow-up time, i.e., the lifespan $\tilde T(\epsilon)$, we
know that
\begin{eqnarray}\phi_y(t,y)>0,\;\;0\le t<\tilde T(\epsilon).\end{eqnarray} The Cauchy problem (1.1) and (1.11) has a unique $C^1$
smooth solution $u=u(t,x)$ and the transformation $(t,y)\rightarrow
(t,x): (t,x)=(t, \phi(t,y))$ is $C^1$ invertible before the lifespan
$\tilde T(\epsilon)$. Therefore, $(\phi,z)=(\phi (t,y),z(t,y))$ is
$C^1$ smooth when the time $0\le t< \tilde T(\epsilon)$. It is
obvious that (4.6) can be deduced from (4.5). Thus, in order to
prove Theorem 4.1, it suffices to prove (4.5) when $0\le t<
\min\left\{\tilde T(\epsilon),\;M_1\epsilon^{-(\alpha+1)}\right\}$.
To do so, it is sufficient to give uniform {\it a priori} estimates
of $C^1$ norm of $z=z(t,y)$ and $\phi=\phi(t,y)$ in the domain
$D(M_1)$.

We fix that
\begin{eqnarray}0<\tau_{1}=\min\left\{\epsilon^{\alpha+1}\tilde T(\epsilon),\;M_1\right\}\le M_1
\end{eqnarray}
and introduce
\begin{eqnarray}k=\partial_y\phi,\;\;\tilde w_i=l_i(z)\partial_yz=w_ik.\end{eqnarray}

Assume that $$z(t,y)=\epsilon\sigma(\tau,y),$$ where we denote
$$\tau=\epsilon^{\alpha+1}t.$$
Introducing the supplemental invariants
\begin{eqnarray}\zeta_i=l_i\left(\epsilon\sigma\right)\partial_y\sigma,\quad\zeta_j=\epsilon^{\alpha+1}l_j\left(\epsilon\sigma\right)\partial_\tau\sigma \quad ( j\neq i),\end{eqnarray}
by (4.3)-(4.4) we have
\begin{eqnarray}\epsilon^{\alpha+1}\partial_\tau\sigma=\sum_{j\neq
i}\zeta_jr_j(\epsilon\sigma)+\epsilon^{-1}
b_i(\epsilon\sigma)r_i(\epsilon\sigma),\;\partial_y\sigma=\sum_{j\neq
i}\frac {
k(\epsilon^{-1}b_j(\epsilon\sigma)-\zeta_j)}{\lambda_j\left(\epsilon\sigma\right)-\lambda_i\left(\epsilon\sigma\right)}r_j\left(\epsilon\sigma\right)+\zeta_ir_i(\epsilon\sigma).\end{eqnarray}
We denote  $\tilde\zeta
=(\zeta_1,\cdots,\widehat{\zeta_i},\cdots,\zeta_n)^t$ and $\tilde b
=(b_1,\cdots,\widehat{b_i},\cdots,b_n)^t$ which do not include
$\zeta_i$ and $b_i$ respectively.

By (1.19) in Theorem 1.1, we have
\begin{eqnarray}|z(t,y)|=|\epsilon\sigma(t,y)|=|u(t,\phi(t,y))|\le U_{\infty}(t)\le
K_3\epsilon,\;\;\mbox{when}\;\;(t,y)\in D(M_1).\end{eqnarray}

We now estimate $k$, $\zeta_i$ and $\tilde \zeta$. Denote
$$K(T)=\;\underset{0\le t\le T\le
\tau_1\epsilon^{-(\alpha+1)}}{\max}\;\;\underset{y\in\mathbb{R}}{\sup}\;|k(t,y)|,$$
$$H^{(i)}(\overline{\tau})=\;\underset{0\le \tau\le \overline{\tau}<
\tau_1\le
M_1}{\max}\;\;\underset{y\in\mathbb{R}}{\sup}\;|\zeta_i(\tau,y)|,$$
$$\tilde H(\overline{\tau})=\;\underset{0\le \tau\le \overline{\tau}<
\tau_1\le M_1}{\max}\;\;\underset{y\in\mathbb{R}}{\sup}\;|\tilde
\zeta(\tau,y)|.$$ It is obvious that
$$K(0)\equiv 1,\quad H^{(i)}(0)=O(1),\quad  \tilde H(0)=O(1).$$

It follows from (4.1) that
\begin{eqnarray}
\partial_tk=\phi_{ty}=(\lambda_i(u))_y=\nabla\lambda_i(u)u_xk=\nabla\lambda_i(u)\left(\sum_{j=1}^nw_jr_j(u)\right)k.
\end{eqnarray} On the other hand, by the Hadamard's formula we have
\begin{eqnarray*}\nabla\lambda_i(u)r_i(u)&=&\left(\nabla\lambda_i(u)r_i(u)-\nabla\lambda_i(u_ie_i)r_i(u_ie_i)\right)+\nabla\lambda_i(u_ie_i)r_i(u_ie_i)
\nonumber\\&=&\sum_{j\ne i}\left[\int_0^1\frac{\partial
(\nabla\lambda_ir_i)}{\partial
u_j}(su_1,\cdots,su_{i-1},u_i,su_{i+1},\cdots,su_n)ds\right]u_j+\nabla\lambda_i(u_ie_i)r_i(u_ie_i).
\end{eqnarray*}
Noting (1.7), (1.15) and (1.19), we obtain
\begin{eqnarray*}
|\nabla\lambda_i(u)r_i(u)|\le C_{20}\left(\sum_{j\ne
i}|v_j|+\epsilon^{\alpha}\right).
\end{eqnarray*}
It is obvious that \begin{eqnarray}w_i(t,\phi(t,y))k(t,y)=\tilde
w_i(t,\phi(t,y))=\epsilon\zeta_i(\tau,
y),\;\;\mbox{where}\;\;\tau=\epsilon^{\alpha+1}t.\end{eqnarray}
Then, it follows from Theorem 1.1 that
\begin{eqnarray*}
K(T)&\le& 1+C_{21}\left[\tilde W_1(T)K(T)+\epsilon
H^{(i)}(\epsilon^{\alpha+1}T)\left(\tilde
V_1(T)+\epsilon^{-1}M_1\right)\right]\\&\le& 1+C_{22}\left[\epsilon
K(T)+\epsilon^2
H^{(i)}(\epsilon^{\alpha+1}T)+M_1H^{(i)}(\epsilon^{\alpha+1}T)\right].
\end{eqnarray*}
Therefore, we get
\begin{eqnarray}
K(T)\le 2+C_{23}M_1H^{(i)}(\epsilon^{\alpha+1}T).
\end{eqnarray}

From (2.7), we have
\begin{eqnarray*}
\frac{d\tilde w_i}{d_it}&=&G_i(t,x)k+w_i\partial_tk(t,y)\\
&=&G_i(t,x)k+w_i\nabla\lambda_i(u)\left(\sum_{j=1}^nw_jr_j(u)\right)k\\
&=&\left[\sum_{l\ne
i}\left(\gamma_{iil}(u)+\gamma_{ili}(u)+\nabla\lambda_i(u)r_l(u)\right)w_l\right.\\&&\left.+\sum_{l=1}^n\sum_{p\ne
q}\sigma_{iil}(u)b_{lpq}(u)u_pu_q+(\tilde b_{ii}(u)-\tilde
b_{ii}(u_ie_i))\right]\tilde w_i\\
&&+\underset{j\ne l;j,l\ne i}{\sum}\gamma_{ijl}(u)w_jw_lk+\sum_{j\ne
i}\sum_{l=1}^n\sigma_{ijl}(u)b_l(u)w_jk+\sum_{l\ne i}(\tilde
b_{il}(u)-\tilde b_{il}(u_le_l))w_lk\\
&\triangleq& a(t,\phi(t,y))\tilde w_i+b(t,\phi(t,y))k(t,y).
\end{eqnarray*}
Thus, we get
\begin{eqnarray*}
\tilde
w_i(t,\phi(t,y))=\left[\int_0^tb(s,\phi(s,y))k(t,y)\exp\left(-\int_0^sa(s^{\prime},\phi(s^{\prime},y))ds^{\prime}\right)ds+
\tilde w_i(0,\phi(0,y))\right]\\
\times\exp\left(\int_0^ta(s^{\prime},\phi(s^{\prime},y))ds^{\prime}\right).
\end{eqnarray*}
From Remark 4.1, (4.12)-(4.13) and (2.22), we know that
\begin{eqnarray}
|w_j(t,\phi(t,y))|\le C_{24}|\partial_t z(t,y)|\le
C_{24}\epsilon^{\alpha+2}|\partial_\tau\sigma(\epsilon^{\alpha+1}t,y)|\le
C_{25}\epsilon(\tilde
H(\epsilon^{\alpha+1}t)+\epsilon),\quad\forall\;j\ne
i.\end{eqnarray} Thus, Theorem 1.1 implies that
\begin{eqnarray*}
\int_0^t|a(s,\phi(s,y))|ds\le C_{26}\left[\tilde
W_1(t)+V_{\infty}(t)\tilde V_1(t)+\tilde V_1(t)\right]\le
C_{27}\epsilon
\end{eqnarray*}
and
\begin{eqnarray*}
\int_0^t|b(s,\phi(s,y))|ds&\le& C_{28}\left[C_{25}\epsilon(\tilde
H(\epsilon^{\alpha+1}t)+\epsilon)\tilde W_1(t)+V_{\infty}(t)^2\tilde
W_1(t)+\tilde W_1(t)V_{\infty}(t)\right]\\&\le&
C_{29}(\epsilon^2+\epsilon^2\tilde H(\epsilon^{\alpha+1}t)).
\end{eqnarray*}
Therefore, noting (4.15) we obtain
\begin{eqnarray}
H^{(i)}(\bar\tau)\le C_{30}\left[H^{(i)}(0)+\epsilon (1+\tilde
H(\bar\tau))K(\epsilon^{-(\alpha+1)}\bar\tau)\right].
\end{eqnarray}

$(4.4)$ gives
\begin{eqnarray*}l_j\left(\epsilon\sigma\right)\left[\epsilon^{\alpha+1}k\partial_\tau
\sigma+\left(\lambda_j\left(\epsilon\sigma\right)-\lambda_i\left(\epsilon\sigma\right)\right)\partial_y\sigma\right]=\epsilon^{-1}b_j(\epsilon\sigma)k,\quad
\forall\; j\neq i.\end{eqnarray*} Differentiating it with respect to
$\tau$ and then multiplying $\epsilon^{\alpha+1}$ yields
\begin{eqnarray*}
&&\epsilon\left(\nabla_ul_j(\epsilon\sigma)(\epsilon^{\alpha+1}\partial_\tau\sigma)\right)^t\left[k(\epsilon^{\alpha+1}\partial_\tau\sigma)
+(\lambda_j(\epsilon\sigma)-\lambda_i(\epsilon\sigma))\partial_y\sigma
\right]\\
&&+l_j(\epsilon\sigma)\left[\epsilon^{\alpha+1}
k\partial_\tau+(\lambda_j(\epsilon\sigma)-\lambda_i(\epsilon\sigma))\partial_y\right](\epsilon^{\alpha+1}\partial_\tau\sigma) \\
&&+ l_j(\epsilon\sigma)\left[(\epsilon^{\alpha+1}\partial_\tau
k)(\epsilon^{\alpha+1}\partial_\tau\sigma)
+\epsilon\left((\nabla_u\lambda_j(\epsilon\sigma)-\nabla_u\lambda_i(\epsilon\sigma))(\epsilon^{\alpha+1}\partial_\tau\sigma)\right)^t\partial_y\sigma\right]\\
&=&\epsilon^{\alpha}\partial_\tau(b_j(\epsilon\sigma)k).
\end{eqnarray*}
Furthermore, we have
\begin{eqnarray*}
&&\left[\epsilon^{\alpha+1}
k\partial_\tau+(\lambda_j(\epsilon\sigma)-\lambda_i(\epsilon\sigma))\partial_y\right]\zeta_j\\
&&-\epsilon\left[
k\nabla_ul_j(\epsilon\sigma)(\epsilon^{\alpha+1}\partial_\tau\sigma)
+(\lambda_j(\epsilon\sigma)-\lambda_i(\epsilon\sigma))(\nabla_ul_j(\epsilon\sigma)\partial_y\sigma)\right]^t(\epsilon^{\alpha+1}\partial_\tau\sigma)\\
&&+\epsilon\left(\nabla_ul_j(\epsilon\sigma)(\epsilon^{\alpha+1}\partial_\tau\sigma)\right)^t\left[k(\epsilon^{\alpha+1}\partial_\tau\sigma)
+(\lambda_j(\epsilon\sigma)-\lambda_i(\epsilon\sigma))\partial_y\sigma
\right]\\
&&+ l_j(\epsilon\sigma)\left[(\epsilon^{\alpha+1}\partial_\tau
k)(\epsilon^{\alpha+1}\partial_\tau\sigma)
+\epsilon\left((\nabla_u\lambda_j(\epsilon\sigma)-\nabla_u\lambda_i(\epsilon\sigma))(\epsilon^{\alpha+1}\partial_\tau\sigma)\right)^t\partial_y\sigma\right]\\
&=&\epsilon^{\alpha}\partial_\tau(b_j(\epsilon\sigma)k).
\end{eqnarray*}
We can rewrite it in a simple form as follows
\begin{eqnarray}
&&\left[\epsilon^{\alpha+1}
k\partial_\tau+(\lambda_j(\epsilon\sigma)-\lambda_i(\epsilon\sigma))\partial_y\right]\zeta_j\nonumber\\
&&+\epsilon\left[(\epsilon^{\alpha+1}\partial_\tau\sigma)^tQ^{(1)}(\epsilon\sigma)(\epsilon^{\alpha+1}\partial_\tau\sigma)k
+(\epsilon^{\alpha+1}\partial_\tau\sigma)^tQ^{(2)}(\epsilon\sigma)(\partial_y\sigma)\right]\nonumber\\
&=&\epsilon^{\alpha}\partial_\tau(b_j(\epsilon\sigma)k)-\epsilon^{\alpha+1}\partial_\tau
k\zeta_j,
\end{eqnarray}
here and hereafter $Q^{(p)}(\epsilon\sigma)\;(p=1,2,\cdots,9)$ are
matrix, column vectors or scalar quantities which are dependent of
$\epsilon\sigma$ continuously.

On the other hand, it follows from (4.12) that
\begin{eqnarray*}
&&(\epsilon^{\alpha+1}\partial_\tau\sigma)^tQ^{(1)}(\epsilon\sigma)(\epsilon^{\alpha+1}\partial_\tau\sigma)k
+(\epsilon^{\alpha+1}\partial_\tau\sigma)^tQ^{(2)}(\epsilon\sigma)(\partial_y\sigma)\\
&=&\left(\sum_{j\ne
i}\zeta_jr_j(\epsilon\sigma)+\epsilon^{-1}b_ir_i(\epsilon\sigma)\right)^tQ^{(1)}(\epsilon\sigma)\left(\sum_{l\ne
i}\zeta_lr_l(\epsilon\sigma)+\epsilon^{-1}b_ir_i(\epsilon\sigma)\right)k\\
&&+\left(\sum_{j\ne
i}\zeta_jr_j(\epsilon\sigma)+\epsilon^{-1}b_ir_i(\epsilon\sigma)\right)^tQ^{(2)}(\epsilon\sigma)\\
&&\times\left(\sum_{l\ne
i}\frac{-\zeta_l}{\lambda_l(\epsilon\sigma)-\lambda_i(\epsilon\sigma)}kr_l(\epsilon\sigma)+\sum_{l\ne
i}\frac{\epsilon^{-1}b_l}{\lambda_l(\epsilon\sigma)-\lambda_i(\epsilon\sigma)}kr_l(\epsilon\sigma)
+\zeta_ir_i(\epsilon\sigma)\right)\\
&=&k\left[\tilde\zeta^t
Q^{(3)}(\epsilon\sigma)\tilde\zeta+\epsilon^{-1}Q^{(4)}(\epsilon\sigma)^t\tilde\zeta
b_i+\epsilon^{-2}Q^{(5)}(\epsilon\sigma)b_i^2+\epsilon^{-2}Q^{(6)}(\epsilon\sigma)^t\tilde
b b_i+\epsilon^{-1}\tilde\zeta^tQ^{(7)}(\epsilon\sigma)\tilde
b\right]\\
&&+\epsilon^{-1}Q^{(8)}(\epsilon\sigma)b_i\zeta_i+Q^{(9)}(\epsilon\sigma)^t\tilde\zeta\zeta_i.
\end{eqnarray*}
Therefore, noting (2.22) and (4.13) we have
$|b_j(\epsilon\sigma)|\le
C_{31}\epsilon^2\;(\forall\;j=1,2,\cdots,n)$ and then from (4.16)
and (4.18) we get, $\forall\;\tau\in [0,\bar\tau]\subset[0,\tau_1)$,
\begin{eqnarray}
&&\left|(\epsilon^{\alpha+1}\partial_\tau\sigma)^tQ^{(1)}(\epsilon\sigma)(\epsilon^{\alpha+1}\partial_\tau\sigma)(\tau,y)k(\epsilon^{-(\alpha+1)}\bar\tau,y)
+(\epsilon^{\alpha+1}\partial_\tau\sigma)^tQ^{(2)}(\epsilon\sigma)(\partial_y\sigma)(\tau,y)\right|\nonumber\\
&\le& C_{32}\left[\left(\tilde H(\bar\tau)^2+\epsilon \tilde
H(\bar\tau)+\epsilon^2\right)K(\epsilon^{-(\alpha+1)}\bar\tau)+\epsilon
H^{(i)}(\bar\tau)+\tilde H(\bar\tau)H^{(i)}(\bar\tau)\right].
\end{eqnarray}
In the proof of (4.16), we have deduced that
\begin{eqnarray}\int_0^{\bar\tau}|\partial_\tau
k(t,y)|d\tau&=&\int_0^{\epsilon^{-(\alpha+1)}\bar\tau}|\partial_tk(t,y)|dt\nonumber\\
&\le&C_{22}\left[\epsilon
K(\epsilon^{-(\alpha+1)}\bar\tau)+\epsilon^2
H^{(i)}(\bar\tau)+M_1H^{(i)}(\bar\tau)\right]\nonumber\\
&\le&C_{33}\left(\epsilon +M_1H^{(i)}(\bar\tau)\right).
\end{eqnarray}
On the other hand, we easily get
\begin{eqnarray}
\left|\int_0^{\bar\tau}(b_j(\epsilon\sigma)k)_\tau
d\tau\right|=\left|\int_0^{\epsilon^{-(\alpha+1)}\bar\tau}(b_j(\epsilon\sigma)k)_t
dt\right|\le C_{34}\epsilon^2K(\epsilon^{-(\alpha+1)}\bar\tau).
\end{eqnarray}

(4.8) ensures that the family of $i-th$ characteristics do not
collapse. Integrating (4.19) and making use of (4.20)-(4.22), we
obtain
\begin{eqnarray}
\tilde H(\bar\tau)
&\le&\tilde H(0)+C_{34}\epsilon^{\alpha+2}K(\epsilon^{-(\alpha+1)}\bar\tau)\nonumber\\
&&+C_{32}\epsilon\left[\left(\tilde H(\bar\tau)^2+\epsilon \tilde
H(\bar\tau)+\epsilon^2\right)K(\epsilon^{-(\alpha+1)}\bar\tau)+\epsilon
H^{(i)}(\bar\tau)+\tilde H(\bar\tau)H^{(i)}(\bar\tau)\right]M_1\nonumber\\
&&+C_{33}\epsilon^{\alpha+1}\tilde H(\bar\tau)\left(\epsilon
+M_1H^{(i)}(\bar\tau)\right).\end{eqnarray}

Therefore, we can deduce from (4.16), (4.18) and (4.23) that,
$\forall\;\;T\in[0,\tau_1\epsilon^{-(\alpha+1)}),\;\;\bar\tau\in[0,\tau_1)$,
\begin{eqnarray}
K(T)\le 2+2C_{23}C_{30}M_1H^{(i)}(0),\;\;H^{(i)}(\bar\tau)\le
2C_{30}H^{(i)}(0),\;\;\tilde H(\bar\tau)\le 2\tilde H(0),
\end{eqnarray}
provided that $\epsilon$ is sufficiently small.

Furthermore, noting (4.12)-(4.13) it follows from (4.24) that
\begin{eqnarray}
|\sigma|\le C_{35},\quad|\partial_\tau\sigma|\le
C_{35}\epsilon^{-(\alpha+1)},\quad|\partial_y\sigma|\le C_{35}
\end{eqnarray}
and
\begin{eqnarray}
|z|\le C_{35}\epsilon,\quad|\partial_tz|\le
C_{35}\epsilon,\quad|\partial_yz|\le C_{35}\epsilon.
\end{eqnarray}
Therefore, (4.6) and (4.7) hold. Then (4.14)-(4.15) and (4.24) imply
that
\begin{eqnarray}
|\phi_y|\le C_{36},\quad|\phi_{ty}|\le
C_{36}\epsilon,\quad|\phi_t|=|\lambda_i(z)|\le C_{36}.
\end{eqnarray}

Therefore, the estimates in (4.5) can be deduced from (4.26)-(4.27).
Thus, Theorem 4.1 is completely proved.
$\quad\quad\quad\quad\quad\quad\quad\quad\quad\Box$

\section{Estimate of lifespan---Proof of Theorem 1.2}

In order to prove Theorem 1.2, i.e., (1.26), as we already assume
that $u$ are the normalized coordinates and noting (1.24)-(1.25), it
suffices to prove
\begin{eqnarray}
\lim_{\epsilon\longrightarrow 0}(\epsilon^{\alpha+1}\tilde
T(\epsilon))=M_0, \end{eqnarray} where
\begin{eqnarray}
M_0&=& \left\{\max_{i\in J_1}\sup_{x\in \mathbb{R}}\left\{-\frac
{1}{\alpha !}\frac{d^ {\alpha+1}\lambda_i}{ds^{\alpha+1}}(0)\psi_i
(x)^{\alpha}\psi^{\prime}_i(x)\right\}\right\}^{-1}\nonumber\\&=&\left\{\max_{i\in
\{1,2,\cdots,n\}}\sup_{x\in \mathbb{R}}\left\{-\frac {1}{\alpha
!}\frac{d^ {\alpha+1}\lambda_i}{ds^{\alpha+1}}(0)\psi_i
(x)^{\alpha}\psi^{\prime}_i(x)\right\}\right\}^{-1}.
\end{eqnarray}
In order to prove (5.1), similar to L. H\"ormander \cite{h} and Kong
\cite{k1}-\cite{k2}, it suffices to show that

\noindent (I). for any fixed $M^*\ge M_0$, we have $\tilde
T(\epsilon)\le M^*\epsilon^{-(\alpha+1)}$, namely
$\underset{\epsilon\longrightarrow 0}{\overline{\lim}}
(\epsilon^{\alpha+1}\tilde T(\epsilon))\le M_0$

and

\noindent (II). for any fixed $M_*\le M_0-\epsilon^{\frac{1}{2}}$,
we have $\tilde T(\epsilon)\ge M_*\epsilon^{-(\alpha+1)}$, namely
$\underset{\epsilon\longrightarrow
0}{\underline{\lim}}(\epsilon^{\alpha+1}\tilde T(\epsilon))\ge M_0$.

Let
\begin{eqnarray}
T_*=K_*\epsilon^{-(\alpha+1)},\;\;T^*=M^*\epsilon^{-(\alpha+1)},
\end{eqnarray}
where $K_*$ is the positive constant $K_6$ given in (1.22) and $M^*$
is an arbitrary fixed constant satisfying that $M^*\ge M_0$. It is
easy to see that
\begin{eqnarray}
0<T_*<T^*<K_4\epsilon^{-(\alpha+\frac{3}{2})}.
\end{eqnarray}

Let $x=x_i(t,y)$ $(i=1,2,\cdots,n)$ be the $i-th$ characteristic
passing through an arbitrary given point $(0,y)$. On any given
existence domain $0\le t\le T\;(T\le T^*)$ of the $C^1$ solution
$u=u(t,x)$, we consider (2.7) along the $i-th$ characteristic
$x=x_i(t,y)$. We can rewrite (2.7) as
\begin{eqnarray}
\frac{dw_i}{d_it}=a_0(t;i,y)w_i^2+a_1(t;i,y)w_i+a_2(t;i,y),
\end{eqnarray}
where
\begin{eqnarray}
a_0(t;i,y)&=&\gamma_{iii}(u),\\
a_1(t;i,y)&=&\sum_{j\ne
i}(\gamma_{iij}(u)+\gamma_{iji}(u))w_j+\sum_{j=1}^n\sigma_{iji}(u)b_j(u)+(\tilde
b_{ii}(u)-\tilde b_{ii}(u_ie_i)),\\
a_2(t;i,y)&=&\sum_{j\ne k;j,k\ne
i}\gamma_{ijk}(u)w_jw_k+\sum_{j=1}^n\sum_{k\ne
i}\sigma_{ijk}(u)b_j(u)w_k+\sum_{k\ne i}\tilde b_{ik}(u)w_k,
\end{eqnarray}
in which $u=u(t,x_i(t,y))$ and
$w_j=w_j(t,x_i(t,y))\;(j=1,2,\cdots,n)$.

\begin{Lemma}
On any given domain $0\le t\le T(\le T^*)$ of the $C^1$ solution
$u=u(t,x)$, there exist positive constants $K_8$ independent of
$\epsilon,y$ and $T$ such that the following estimates hold:
\begin{eqnarray}
\int_0^T|a_1(t;i,y)|dt&\le& K_8\epsilon,\\
\int_0^T|a_2(t;i,y)|dt&\le& K_8\epsilon^2,\\
K(i,y;0,T)&\triangleq&\int_0^T|a_2(t;i,y)|dt\cdot\exp
\left(\int_0^T|a_1(t;i,y)|dt\right)\le K_8\epsilon^2.
\end{eqnarray}
\end{Lemma}
\noindent{{\bf Proof.}} It follows from (2.18)-(2.24) and Theorem
1.1 that
\begin{eqnarray}
\int_0^T|a_1(t;i,y)|dt&\le& C_{37}(\tilde W_1(T)+V_{\infty}(T)\tilde
V_1(T)+\tilde V_1(T))\le C_{38}\epsilon.
\end{eqnarray}

In Theorem 4.1, we take $M_1=M^*+1$. Noting (4.7) and (1.19), we
easily see that
\begin{eqnarray*}\sum_{j\ne k;j,k\ne
i}\int_0^T|w_jw_k|(s,x_i(s,y))ds&\le& C_{19}(n-2)\epsilon\sum_{k\ne
i}\int_0^T|w_k|(s,x_i(s,y))ds \\&\le& C_{39}(n-2)\epsilon\sum_{k\ne
i}\tilde W_1(T)\;\\&\le&\; C_{40}\epsilon^2.\end{eqnarray*}

Therefore, noting (2.22), (2.24) and Theorem 1.1, we get
\begin{eqnarray}
\int_0^T|a_2(t;i,y)|dt&=& \int_0^T\left|\sum_{j\ne k;j,k\ne
i}\gamma_{ijk}(u)w_jw_k+\sum_{j=1}^n\sum_{k\ne
i}\sigma_{ijk}(u)\sum_{p\ne
q}b_{jpq}(u)u_pu_qw_k\right.\nonumber\\&&\left.\;\;\;\;\;+\sum_{k\ne
i}(\tilde b_{ik}(u)-\tilde
b_{ik}(u_ke_k))w_k\right|dt\nonumber\\
&\le&C_{41}\left[\sum_{j\ne k;j,k\ne
i}\int_0^T|w_jw_k|(s,x_i(s,y))ds+U_{\infty}(T)^2\tilde
W_1(T)+U_{\infty}(T)\tilde
W_1(T)\right]\nonumber\\
&\le&C_{42}\epsilon^2.
\end{eqnarray}
Then (5.9)-(5.11) can be easily deduced from
(5.12)-(5.13).$\quad\Box$

Similar to Lemma 1.4.1 in L. H\"ormander \cite{h}, we have
\begin{Lemma}
Let $z=z(t)$ be a solution in $[0,T]$ of the Riccadi's differential
equation:
\begin{eqnarray*}
\frac{dz}{dt}=a_0(t)z^2+a_1(t)z+a_2(t),
\end{eqnarray*}
where $a_j(t)\;(j=0,1,2)$ are continuous and $T>0$ is a given real
number. Let \begin{eqnarray*} K=\int_0^T|a_2(t)|dt\cdot
\exp\left(\int_0^T|a_1(t)|dt\right).\end{eqnarray*} If $z(0)>K$,
then it follows that
\begin{eqnarray*} \int_0^T|a_0(t)|dt\cdot\exp\left(-\int_0^T|a_1(t)|dt\right)<(z(0)-K)^{-1}.\end{eqnarray*}
\end{Lemma}
\begin{Remark} L. H\"ormander assumed that $a_0(t)\ge 0$ in Lemma
1.4.1 (see page 230 in \cite{h}). In Lemma 5.2, we do not assume
this. It is easy to find that we can prove Lemma 5.2 similar to
Lemma 1.4.1 in L. H\"ormander \cite{h}.
\end{Remark}

Next we give the estimate of the lifespan of classical solution to
the Cauchy problem (1.1) and (1.11) under the assumptions of Theorem
1.2.

\noindent{{\bf (I) Upper bound of the lifespan--- Estimate on
$\underset{\epsilon\longrightarrow
0^+}{\overline{\lim}}(\epsilon^{\alpha+1}\tilde T(\epsilon))\le
M_0$}}

It follows from (2.6), (2.18)-(2.24) and (1.19) that, along the
$i-th$ characteristic $x=x_i(t,y)$,
\begin{eqnarray*}
|v_i(i;t,y)-v_i(i;0,y)|&\le&\int_0^t|F_i(s,x_i(s,y))|ds\\
&\le & \int_0^t\left|\sum_{k\ne
i}\sum_{j=1}^n\beta_{ijk}(u)v_jw_k+\sum_{j,k}\nu_{ijk}(u)v_j\sum_{p\ne
q}b_{kpq}(u)u_pu_q\right.\\&&\;\;\;\;\left.+\sum_{p\ne q}b_{ipq}(u)u_pu_q\right|(s,x_i(s,y))ds\\
&\le& C_{43}\left[V_{\infty}(t)\tilde W_1(t)+V_{\infty}(t)^2\tilde
V_1(t)+V_{\infty}(t)\tilde V_1(t)\right]\\
&\le&C_{44}\epsilon^2.
\end{eqnarray*}
Then, as $u$ are the normalized coordinates and $l_i(0)=e_i$, from
(1.19) we easily get, along the $i-th$ characteristic $x=x_i(t,y)$,
\begin{eqnarray*}
|u_i(i;t,y)-u_i(i;0,y)|=|u_i(i;t,y)-f_i(\epsilon,y)|\le
C_{45}\epsilon^2.
\end{eqnarray*}

Using Hadamard's formula and noting (1.12)-(1.13), from (5.6) we
get, along the $i-th$ characteristic $x=x_i(t,y)$,
\begin{eqnarray}
a_0(t;i,y)&=&\gamma_{iii}(u)=\gamma_{iii}(u_ie_i)+(\gamma_{iii}(u)-\gamma_{iii}(u_ie_i))\nonumber \\
&=&\gamma_{iii}(u_ie_i)+\sum_{j\ne i}\left[\int_0^1\frac{\partial
\gamma_{iii}}{\partial
u_j}(su_1,\cdots,su_{i-1},u_i,su_{i+1},\cdots,su_n)ds\right]u_j\nonumber\\
&=&-\frac{1}{\alpha !}\frac{\partial^{1+\alpha}\lambda_i}{\partial
u_i^{1+\alpha}}(0)(u_i)^{\alpha}+O(\epsilon^{1+\alpha})\nonumber\\
&&+\sum_{j\ne i}\left[\int_0^1\frac{\partial \gamma_{iii}}{\partial
u_j}(su_1,\cdots,su_{i-1},u_i,su_{i+1},\cdots,su_n)ds\right]u_j\nonumber\\
&=&-\frac{1}{\alpha !}\frac{\partial^{1+\alpha}\lambda_i}{\partial
u_i^{1+\alpha}}(0)(\epsilon\psi_i(y))^{\alpha}+\alpha O(\epsilon^{\alpha+r})+O(\epsilon^{1+\alpha})\nonumber\\
&&+\sum_{j\ne i}\left[\int_0^1\frac{\partial \gamma_{iii}}{\partial
u_j}(su_1,\cdots,su_{i-1},u_i,su_{i+1},\cdots,su_n)ds\right]u_j
\end{eqnarray}

Noting that the initial data satisfies (1.25), we observe that there
exist an index $i_0\in J_1$ and a point $x_0\in \mathbb{R}$ such
that
\begin{eqnarray}
M_0&=& \left\{-\frac {1}{\alpha !}\frac{\partial^
{\alpha+1}\lambda_{i_0}}{\partial u_{i_0}^{\alpha+1}}(0)\psi_{i_0}
(x_0)^{\alpha}\psi^{\prime}_{i_0}(x_0)\right\}^{-1}.
\end{eqnarray}
Noting (2.10) and (1.15), we have
\begin{eqnarray*}
\frac{\partial ^l\gamma_{i_0i_0i_0}}{\partial
u_{i_0}^l}(0)=0\;\;(l=0,1,\cdots,\alpha-1)\;\;\mbox{but}\;\;\frac{\partial
^{\alpha}\gamma_{i_0i_0i_0}}{\partial u_{i_0}^{\alpha}}(0)\ne 0.
\end{eqnarray*}
Then (5.15) becomes
\begin{eqnarray}
M_0&=& \left\{\frac {1}{\alpha !}\frac{\partial^
{\alpha}\gamma_{i_0i_0i_0}}{\partial u_{i_0}^{\alpha}}(0)\psi_{i_0}
(x_0)^{\alpha}\psi^{\prime}_{i_0}(x_0)\right\}^{-1}\triangleq
(b\psi_{i_0}^{\prime}(x_0))^{-1}.
\end{eqnarray}
Without loss of generality, we may suppose that
\begin{eqnarray}
b>0\;\; and\;\;\psi_{i_0}^{\prime}(x_0)>0.
\end{eqnarray}
Otherwise, changing the sign of $u$, we can draw the same
conclusion.

Noting (1.12)-(1.13), (1.25) and (5.11), we get immediately
\begin{eqnarray*}
w_{i_0}(0,x_0)&=&l_{i_0}\left(f(\epsilon,x_0)\right)\frac{\partial
f}{\partial x}(\epsilon,x_0)\\
&=&\left[l_{i_0}(0)+O(\epsilon)\right]\times\left[\frac{\partial
f}{\partial x}(0,x_0)+\frac{\partial^2f}{\partial \epsilon\partial
x}(0,x_0)\epsilon+O(\epsilon^{1+r})\right]\\
&=&\epsilon\psi_{i_0}^{\prime}(x_0)+O(\epsilon^{1+r})>K_8\epsilon^2\ge
K(i_0,x_0;0,T).
\end{eqnarray*}
Therefore, we immediately observe that Lemma 5.2 (revised version of
Lemma 1.4.1 in L. H\"ormander \cite{h}) can be applied to the
initial value problem for (5.5) with the following initial condition
\begin{eqnarray}
t=0:\;\;w_{i_0}=w_{i_0}(0,x_0)=\epsilon\psi_{i_0}^{\prime}(x_0)+O(\epsilon^{1+r})
\end{eqnarray}
and then we obtain
\begin{eqnarray*}
\int_0^T|a_0(t;i_0,x_0)|dt\cdot\exp\left(-\int_0^T|a_1(t;i_0,x_0)|dt\right)<\left(w_{i_0}(0,x_0)-K(i_0,x_0;0,T)\right)^{-1},
\end{eqnarray*}
namely,
\begin{eqnarray}
&&\exp\left(-\int_0^T|a_1(t;i_0,x_0)|dt\right)\times\nonumber\\
&&\int_0^T|a_0(t;i_0,x_0)|\left(w_{i_0}(0,x_0)-K(i_0,x_0;0,T)\right)dt<1.
\end{eqnarray}
Substituting (5.14) into (5.19) and noting (1.19) and
the fact that
$T\le T^*=M^*\epsilon^{-(1+\alpha)}$, we obtain
\begin{eqnarray}
\underset{\epsilon\longrightarrow 0}{\overline{\lim}}
\left\{\epsilon^{\alpha+1}T\cdot\frac{1}{\alpha !}\frac{\partial^
{\alpha}\gamma_{i_0i_0i_0}}{\partial u_{i_0}^{\alpha}}(0)\psi_{i_0}
(x_0)^{\alpha}\psi^{\prime}_{i_0}(x_0)\right\}\le 1.
\end{eqnarray}
Noting (5.16), from (5.20) we get immediately
\begin{eqnarray}
\underset{\epsilon\longrightarrow
0}{\overline{\lim}}(\epsilon^{\alpha+1}\tilde T(\epsilon))\le M_0.
\end{eqnarray}
(5.21) gives an upper bound of the lifespan $\tilde T(\epsilon)$.

\noindent{{\bf (II) Lower bound of the lifespan--- Estimate on
$\underset{\epsilon\longrightarrow
0^+}{\underline{\lim}}(\epsilon^{\alpha+1}\tilde T(\epsilon))\ge
M_0$}}

To do so, it suffices to prove that, for any fixed $M_*$ satisfying
that \begin{eqnarray} 0<M_*<M_0-\epsilon^{\frac{1}{2}r},
\end{eqnarray}
we have
\begin{eqnarray}
\tilde T(\epsilon)\ge M_*\epsilon^{-(\alpha+1)},
\end{eqnarray}
provided that $\epsilon>0$ is small enough. Hence, we only need to
establish a uniform {\it a priori} estimate on $C^1$ norm of the
$C^1$ solution $u=u(t,x)$ on any given existence domain $0\le t\le
T\le M_*\epsilon^{-(\alpha+1)}$. The uniform {\it a priori} estimate
on $C^0$ norm of $u=u(t,x)$ has been established in Theorem 1.1. It
remain to establish a uniform {\it a priori} estimate on $C^0$ norm
of the first derivatives of $u=u(t,x)$, namely a uniform {\it a
priori} estimate on $C^0$ norm of
$w=(w_1(t,x),w_2(t,x),\cdots,w_n(t,x))^T$.

In order to estimate $w_i=w_i(t,x)$ on the existence domain $0\le
t\le T$ (where $T$ satisfies $T\le M_*\epsilon^{-(\alpha+1)}$) of
the $C^1$ solution $u=u(t,x)$, we still consider (5.5) along the
$i-th$ characteristic $x=x_i(t,y)$ passing through an arbitrary
fixed point $(0,y)$. Without loss of generality, we may suppose that
\begin{eqnarray}
\psi_i^{\prime}(y)\ge 0.
\end{eqnarray}
Otherwise, changing the sign of $u$, we can draw the same
conclusion.

Let $$a_0^+(t;i,y)=\max\{a_0(t;i,y),0\}.$$ Noting the fact that
$T\le M_*\epsilon^{-(\alpha+1)}$ and using Theorem 1.1, (5.14) and
(5.18), we obtain
\begin{eqnarray}
&&w_i(0,y)\int_0^Ta_0^+(t;i,y)dt\nonumber\\
&\le
&\left(\epsilon\psi_i^{\prime}(y)+O(\epsilon^{1+r})\right)\left\{\max\left\{-\frac{1}{\alpha
!}\frac{\partial^{1+\alpha}\lambda_i}{\partial
u_i^{1+\alpha}}(0)(\epsilon\psi_i(y))^{\alpha},0\right\}T\right.\nonumber\\
&&\left.+C_{46}\left(\alpha
\epsilon^{\alpha+r}T+\epsilon^{1+\alpha}T+\tilde
V_1(T)\right)\right\}\nonumber\\
&\le &\max\left\{-\frac{1}{\alpha
!}\frac{\partial^{1+\alpha}\lambda_i}{\partial
u_i^{1+\alpha}}(0)(\psi_i(y))^{\alpha},0\right\}\psi_i^{\prime}(y)M_*+C_{47}(\epsilon^r+\epsilon)\nonumber\\
&\le &
M_0^{-1}M_*+C_{48}\epsilon^r=M_0^{-1}(M_0-\epsilon^{\frac{1}{2}r})+C_{48}\epsilon^r
< \;1,
\end{eqnarray}
provided that $\epsilon>0$ is small enough. On the other hand,
noting (5.14) and Theorem 1.1, we get immediately
\begin{eqnarray}
\int_0^T|a_0(t;i,y)|dt\le
C_{49}(\epsilon^{\alpha}T+\alpha\epsilon^{\alpha+r}T+\epsilon^{\alpha+1}T+\epsilon)\le
C_{50}M_*\epsilon^{-1}\le C_{51}\epsilon^{-1}.
\end{eqnarray}
Then, noting (5.25)-(5.26) and Lemma 5.1, we obtain
\begin{eqnarray}
\int_0^Ta_0^+(t;i,y)dt\times
\exp\left(\int_0^T|a_1(t;i,y)|dt\right)<(w_i(0,y)+K(i,y;0,T))^{-1}
\end{eqnarray}
and
\begin{eqnarray}
\int_0^T|a_0(t;i,y)|dt\times
\exp\left(\int_0^T|a_1(t;i,y)|dt\right)<(K(i,y;0,T))^{-1},
\end{eqnarray}
where $T\le M_*\epsilon^{-(\alpha+1)}$.

Noting (5.24) and (5.27)-(5.28), we observe that Lemma 1.4.2 in L.
H\"ormander \cite{h} can be applied to the initial value problem for
equation (5.5) with the following initial condition
\begin{eqnarray*}
t=0: w_i=w_i(0,y).
\end{eqnarray*}
Then we have
\begin{eqnarray}
(w_i(T,x_i(T,y)))^{-1}&\ge& (w_i(0,y)+K(i,y;0,T))^{-1}-\int_0^Ta_0^+(t;i,y)dt\nonumber\\
&&\times
\exp\left(\int_0^T|a_1(t;i,y)|dt\right),\;\mbox{if}\;\;w_i(T,x_i(T,y))>0
\end{eqnarray}
and
\begin{eqnarray}
|w_i(T,x_i(T,y))|^{-1}&\ge& (K(i,y;0,T))^{-1}-\int_0^T|a_0(t;i,y)|dt\nonumber\\
&& \times
\exp\left(\int_0^T|a_1(t;i,y)|dt\right),\;\mbox{if}\;\;w_i(T,x_i(T,y))<0.
\end{eqnarray}
Noting (5.25)-(5.26) and Lemma 5.1, from (5.29)-(5.30) we get
respectively
\begin{eqnarray}
(w_i(T,x_i(T,y)))^{-1}\ge
\frac{1}{2}\left(1-\frac{M_*}{M_0}\right)(w_i(0,y)+K(i,y;0,T))^{-1},\;\mbox{if}\;\;w_i(T,x_i(T,y))>0
\end{eqnarray}
and
\begin{eqnarray}
|w_i(T,x_i(T,y))|^{-1} \ge
\frac{1}{2}(K(i,y;0,T))^{-1},\;\mbox{if}\;\;w_i(T,x_i(T,y))<0.
\end{eqnarray}
Therefore, we have
\begin{eqnarray}
w_i(T,x_i(T,y))\le \frac{2}{
1-\frac{M_*}{M_0}}(w_i(0,y)+K(i,y;0,T))\le
C_{52}\epsilon^{1-\frac{1}{2}r},\;\mbox{if}\;\;w_i(T,x_i(T,y))>0
\end{eqnarray}
and
\begin{eqnarray}
|w_i(T,x_i(T,y))|\le 2K(i,y;0,T)\le
2K_8\epsilon^2,\;\mbox{if}\;\;w_i(T,x_i(T,y))<0.
\end{eqnarray}
It follows from (5.33)-(5.34) that
\begin{eqnarray}|w_i(T,x_i(T,y))|\le C_{53}\epsilon^{1-\frac{1}{2}r}.
\end{eqnarray}
For each $i\in\{1,2,\cdots,n\}$ and any $t\in [0,T]$, we can prove
similarly that $w_i(t,x_i(t,y))$ satisfies the same estimate. Noting
that $(0,y)$ is arbitrary, we have
$$||w(t,x)||_{C^0[0,T]\times \mathbb{R}}\le
C_{54}\epsilon^{1-\frac{1}{2}r},$$ where $T\le
M_*\epsilon^{-(\alpha+1)}$. Hence, (5.23) holds and then
\begin{eqnarray}\underset{\epsilon\longrightarrow
0^+}{\underline{\lim}}(\epsilon^{\alpha+1}\tilde T(\epsilon))\ge
M_0.\end{eqnarray} The combination of (5.21) and (5.36) gives
(1.26). Thus, Theorem 1.2 is proved completely. $\quad\Box$

\end{document}